\documentclass[11pt]{article}
\usepackage{amsfonts,amsmath}

\newtheorem{thm}{Theorem}

\newtheorem{lemma}[thm]{Lemma}

\newtheorem{claim}[thm]{Claim}
\usepackage{epsfig}

\def\zed{{\mathbb Z}}

\begin{document}

\def\G{{\Gamma}}
 \def\d{{\delta}}
 \def\ci{{\circ}}
 \def\e{{\epsilon}}
 \def\l{{\lambda}}
 \def\L{{\Lambda}}
 \def\m{{\mu}}
 \def\n{{\nu}}
 \def\o{{\omega}}
  \def\O{{\Omega}}
   \def\Th{{\Theta}}\def\s{{\sigma}}
 \def\v{{\varphi}}
 \def\a{{\alpha}}
 \def\b{{\beta}}
 \def\p{{\partial}}
 \def\r{{\rho}}
 \def\to{{\rightarrow}}
 \def\lra{{\longrightarrow}}
 \def\g{{\gamma}}
 \def\D{{\Delta}}
 \def\La{{\Leftarrow}}
 \def\Ra{{\Rightarrow}}
 \def\x{{\xi}}
 \def\c{{\mathbb C}}
 \def\z{{\mathbb Z}}
 \def\2{{\mathbb Z_2}}
 \def\q{{\mathbb Q}}
 \def\t{{\tau}}
 \def\u{{\upsilon}}
 \def\th{{\theta}}
 \def\Th{{\Theta}}
 \def\la{{\leftarrow}}
 \def\lla{{\longleftarrow}}
 \def\da{{\downarrow}}
 \def\ua{{\uparrow}}
 \def\nwa{{\nwtarrow}}
 \def\swa{{\swarrow}}
 \def\nea{{\netarrow}}
 \def\sea{{\searrow}}
 \def\hla{{\hookleftarrow}}
 \def\hra{{\hookrightarrow}}
 \def\sl2{{SL(2,\mathbb C)}}
 \def\ps{{PSL(2,\mathbb C)}}
 \def\qed{{\hspace{2mm}{\small $\diamondsuit$}}}

 \def\pf{{\noindent{\bf Proof.\hspace{2mm}}}}
 \def\ni{{\noindent}}
 \def\sm{{{\mbox{\tiny M}}}}
  \def\sn{{{\mbox{\tiny N}}}}
\def\smn{{{\mbox{\tiny $M_n$}}}}
 \def\sf{{{\mbox{\tiny F}}}}
  \def\sc{{{\mbox{\tiny C}}}}
 \def\ke{{\mbox{ker}(H_1(\p M;\2)\to H_1(M;\2))}}
 \def\et{{\mbox{\hspace{1.5mm}}}}
 \def\sk{{{\mbox{\tiny K}}}}
 \def\sw{{{\mbox{\tiny W}}}}
 \def\sb{{{\mbox{\tiny \cal B}}}}
 \def\sl{{{\mbox{\tiny $\L$}}}}
 \def\sh{{{\mbox{\tiny $\mathbb H$}}}}
\def\se{{{\mbox{\tiny $(\e_1,\e_2)$}}}}
\def\sg{{{\mbox{\tiny $\g$}}}}
\def\sd{{{\mbox{\tiny $\d$}}}}
\def\st{{{\mbox{\tiny $(2)$}}}}
\def\sa{{{\mbox{\tiny $(W, L, {\cal B})$}}}}

\title{Reducible And Finite Dehn Fillings}

\author{Steven Boyer\footnote{{Partially supported by NSERC grant RGPIN 9446}},
Cameron McA. Gordon, and Xingru Zhang }

\maketitle

\begin{abstract}
We show that the distance between a finite filling slope and a
reducible filling slope on the boundary of a hyperbolic knot
manifold is at most one.
\end{abstract}

Let $M$ be a knot manifold, i.e.  a connected, compact, orientable
$3$-manifold whose boundary is a torus. A knot manifold is said to
be {\it hyperbolic} if its interior admits a complete hyperbolic
metric of finite volume. Let $M(\alpha)$ denote the manifold
obtained by Dehn filling $M$ with slope $\alpha$ and let
$\Delta(\alpha,\beta)$ denote the distance between two slopes
$\alpha$ and $\beta$ on $\partial M$. When $M$ is hyperbolic but
$M(\alpha)$ isn't, we call the corresponding  filling (slope) an
{\it exceptional}  filling (slope). Perelman's recent proof of
Thurston's geometrisation conjecture implies that a filling is
exceptional if and only if it is either reducible, toroidal, or
Seifert fibred. These include all manifolds whose fundamental groups
are either cyclic, finite, or very small (i.e. contain no
non-abelian free subgroup). Sharp upper bounds on the distance
between exceptional filling slopes of various types have been
established in many cases, including: \vspace{-.25cm}
\begin{description}
\vspace{-.15cm} \item[\hspace{5mm} $\bullet$] $\Delta(\alpha,\beta)\leq 1$ if both $\alpha$ and $\beta$ are reducible filling slopes \cite{GL2}
\vspace{-.15cm} \item[\hspace{5mm} $\bullet$] $\Delta(\alpha,\beta)\leq 1$ if both $\alpha$ and $\beta$ are cyclic filling slopes \cite{CGLS}
\vspace{-.15cm} \item[\hspace{5mm} $\bullet$] $\Delta(\alpha,\beta)\leq 1$ if  $\alpha$ is a cyclic filling slope and  $\beta$ is a reducible filling slope \cite{BZ2}
\vspace{-.15cm} \item[\hspace{5mm} $\bullet$] $\Delta(\alpha,\beta)\leq 2$ if  $\alpha$ is a cyclic filling slope and  $\beta$ is a finite filling slope
\cite{BZ1}
\vspace{-.15cm} \item[\hspace{5mm} $\bullet$] $\Delta(\alpha,\beta)\leq 2$ if  $\alpha$ is a reducible filling slope and $\beta$ is a very small filling slope \cite{BCSZ2}
\vspace{-.65cm} \item[\hspace{5mm} $\bullet$] $\Delta(\alpha,\beta)\leq 3$ if both $\alpha$ and $\beta$ are finite filling slopes \cite{BZ3}
\vspace{-.15cm} \item[\hspace{5mm} $\bullet$] $\Delta(\alpha,\beta)\leq 3$ if  $\alpha$ is a reducible filling slope and $\beta$ is a toroidal  filling slope \cite{Wu} \cite{Oh}
\vspace{-.15cm} \item[\hspace{5mm} $\bullet$] $\Delta(\alpha,\beta)\leq 8$ if both $\alpha$  and  $\beta$ are toroidal  filling slopes \cite{Go}
\end{description}
\vspace{-.25cm}
\noindent In this paper we give the sharp upper bound on the distance between
a reducible filling slope and finite filling slope.

\begin{thm}\label{rf}
Let $M$ be a hyperbolic knot manifold.
If $M(\alpha)$ has a finite fundamental group and $M(\beta)$
is a reducible manifold, then $\Delta(\alpha,\beta)\leq 1$.
\end{thm}
Example 7.8 of \cite{BZ2} describes a hyperbolic knot manifold $M$ and
slopes $\alpha_1, \alpha_2, \beta$ on $\partial M$ such that $M(\beta)$ is reducible,
 $\pi_1(M(\alpha_1))$ is finite cyclic, $\pi_1(M(\alpha_2))$ is finite non-cyclic,
 and $\Delta(\alpha_1, \beta) = \Delta(\alpha_2, \beta) = 1$.
In fact there are hyperbolic knot manifolds with reducible and
finite fillings for every finite type: cyclic, dihedral,
tetrahedral, octahedral and icosahedral, in the terminology of
[BZ1]; see [K].

A significant reduction of Theorem \ref{rf} was obtained in \cite{BCSZ2}. Before describing this work, we need to introduce some notation and terminology.

Denote the octahedral group by $O$, the binary octahedral group by $O^*$, and let $\varphi: O^* \to O$ be the usual surjection. We say that $\alpha$ is an {\it $O(k)$-type} filling slope if $\pi_1(M(\alpha)) \cong O^* \times \mathbb Z/j$ for some integer $j$ coprime to $6$ and the image of $\pi_1(\partial M)$ under the composition $\pi_1(M) \to \pi_1(M(\alpha)) \stackrel{\cong}{\longrightarrow} O^* \times \mathbb Z/j \stackrel{{\rm proj}}{\longrightarrow} O^* \stackrel{\varphi}{\longrightarrow} O$ is $\mathbb Z / k$. Clearly $k \in \{1,2,3,4\}$. It is shown in \S 3 of \cite{BZ3} that $k$ is independent of the choice of isomorphism $\pi_1(M(\alpha)) \stackrel{\cong}{\longrightarrow} O^* \times \mathbb Z/j$.

A lens space whose fundamental group has order $p \geq 2$ will be denoted by $L_p$.

\begin{thm} \label{almost}
Let $M$ be a hyperbolic knot manifold.
If $M(\alpha)$ has a finite fundamental group and $M(\beta)$ is a reducible manifold, then $\Delta(\alpha,\beta)\leq 2$. Further, if $\Delta(\alpha,\beta) = 2$, then $H_1(M) \cong \mathbb Z \oplus \mathbb Z/2$, $M(\beta) \cong L_2 \# L_3$, and $\alpha$ is an $O(k)$-type filling slope for some $k \in \{1,2,3\}$.
\end{thm}

\pf This is Theorem 1.1 of \cite{BCSZ2} except that that theorem only claimed that $\alpha$ is an $O(k)$-type filling slope for some $k \in \{1,2,3,4\}$. Since $H_1(M)$ contains $2$-torsion, the argument of the last paragraph of the proof of Theorem 2.3 of \cite{BZ1} (see page 1026 of that article) shows that $k \in \{1,2,3 \}$.
\qed

\noindent Thus, in order to prove Theorem \ref{rf}, we are reduced to considering the case where $H_1(M) \cong \mathbb Z \oplus \mathbb Z/2$, $M(\beta) \cong L_2 \# L_3$, and $\alpha$ is an $O(k)$ type filling slope for some $k \in \{1,2,3 \}$. We do this below. We also assume that $\Delta(\alpha, \beta) = 2$ in order to derive a contradiction.

An {\it essential surface} in $M$ is a compact, connected, orientable, incompressible, and non-boundary parallel, properly embedded $2$-submanifold of $M$. A slope $\beta$ on $\partial M$ is called a {\it boundary slope} if there is an essential surface $F$ in $M$ with non-empty boundary of the given  slope $\beta$. A boundary slope $\beta$ is called {\it strict} if
there is an essential surface $F$ in $M$ of boundary slope $\beta$ such that $F$ is neither  a fiber
nor a semi-fiber. When $M$ has a closed essential surface $S$, let $\mathcal{C}(S)$
be the set of  slopes $\gamma$ on $\partial M$ such that $S$ is compressible in $M(\gamma)$.
A slope $\eta$ is called a {\it singular slope} for $S$ if $\eta \in \mathcal{C}(S)$ and
$\Delta(\eta, \gamma)\leq 1$ for each $\gamma \in \mathcal{C}(S)$.

Since $\pi_1(M(\alpha))$ is finite, the first Betti number of $M$ is $1$, $M(\alpha)$ is irreducible by \cite{GL2}, and neither $\alpha$ nor $\beta$ is a singular slope by  Theorem 1.5 of \cite{BGZ}. As  $M(\beta)$ is reducible, $\beta$ is a boundary slope. Further, by Proposition 3.3 of \cite{BCSZ2} we may assume that up to isotopy, there is a unique
essential surface $P$ in $M$  with boundary slope $\beta$. This surface is necessarily planar. It is also separating as $M(\beta)$ is a rational homology $3$-sphere, and so has an even number of boundary components. This number is at least $4$ since $M$ is hyperbolic.

\begin{lemma} \label{O2}
If $\Delta(\alpha, \beta) = 2$, then $\alpha$ is of type $O(2)$.
\end{lemma}

\pf According to Theorem \ref{rf}, we must show that $\alpha$ does not have type $O(k)$ for $k = 1, 3$.

Let $X_0 \subset X(M(\beta))) \subset X(M)$ be the unique non-trivial curve. (We refer the reader to \S 6 of \cite{BCSZ2} for notation, background results, and further references on $PSL_2(\mathbb C)$ character varieties.) Since $\beta$ is not a singular slope, Proposition 4.10 of \cite{BZ2} implies that the regular function $f_\alpha: X_0 \to \mathbb C, \chi_\rho \mapsto (\hbox{trace}(\rho(\alpha)))^2 - 4$, has a pole at each ideal point of $X_0$. (We have identified $\alpha \in H_1(\partial M)$ with its image in $\pi_1(\partial M) \subset \pi_1(M)$
under the Hurewicz homomorphism.) In particular, the Culler-Shalen seminorm $\|\cdot\|_{X_0}: H_1(\partial M; \mathbb R) \to [0, \infty)$ is non-zero. Hence there is a non-zero integer $s_0$ such that for all $\gamma \in H_1(\partial M)$ we have
$$\|\gamma\|_{X_0} = |\gamma \cdot \beta| s_0$$
where $\gamma \cdot \beta$ is the algebraic intersection number of the two classes (c.f. Identity 6.1.2 of \cite{BCSZ2}). Fix a class $\beta^* \in H_1(\partial M)$ satisfying $\beta \cdot \beta^* = \pm 1$,
so in particular $\|\beta^*\|_{X_0} = s_0$. We can always find such a $\beta^*$ so that
$$\alpha = \beta + 2 \beta^*$$
According to Proposition 8.1 of \cite{BCSZ2}, if $\pm \beta \ne \gamma \in H_1(\partial M)$ is a slope satisfying $\Delta(\alpha, \gamma) \equiv 0$ (mod $k$), then $2s_0 = \|\alpha\|_{X_0} \leq \|\gamma\|_{X_0} = \Delta(\gamma, \beta) s_0$. Hence $\Delta(\gamma, \beta) \geq 2$. Consideration of $\gamma = \beta^*$ and $\gamma = \beta - \beta^*$ then shows that $k \ne 1, 3$.
\qed

\begin{lemma}
If $\Delta(\alpha, \beta) = 2$, then $P$ has exactly four boundary components.
\end{lemma}

\pf We continue to use the notation developed in the proof of the previous lemma.

By Case 1 of \S 8 of \cite{BCSZ2} we have $2 \leq 1 + \frac{3}{s_0} < 3$ and so $s_0$ is either $2$ or $3$. We claim that $s_0 = 2$. To prove this, we shall suppose that $s_0 = 3$ and derive a contradiction.

It follows from the method of proof of Lemma 5.6 of \cite{BZ1} that $\pi_1(M(\alpha)) \cong O^* \times \mathbb Z/j$ has exactly two irreducible characters with values in $PSL_2(\mathbb C)$ corresponding to a representation $\rho_1$ with image $O$ and a representation $\rho_2$ with image $D_3$ (the dihedral group of order 6).
Further, $\rho_2$ is the composition of $\rho_1$ with the quotient of $O$ by
its unique normal subgroup isomorphic to $\mathbb Z/2 \oplus \mathbb Z/2$.
It follows from Proposition 7.6 of \cite{BCSZ2} that if $s_0 = 3$, the characters of $\rho_1$ and $\rho_2$ lie on $X_0$ and provide jumps in the multiplicity of zero of $f_\alpha$ over $f_{\beta^*}$. Lemma 4.1 of \cite{BZ1} then implies that both $\rho_1(\beta^*)$ and $\rho_2(\beta^*)$ are non-trivial. By the previous lemma, $\alpha$ is a slope of type $O(2)$. Thus $\rho_1(\beta^*)$ has order $2$. Since $\rho_2(\beta^*) \ne \pm I$ and $\rho_2$ factors through $\rho_1$, $\rho_2(\beta^*)$ also has order 2.

Next we claim that $\beta^*$ lies in the kernel of the composition of $\rho_2$ with the abelianisation $D_3 \to H_1(D_3; \mathbb Z/2)$. To see this, note first that $\beta$ is non-zero in $H_1(M;\mathbb Z/2) = \mathbb Z/2 \oplus \mathbb Z/2$ since $H_1(M(\beta); \mathbb Z/2) = H_1(L_2 \# L_3; \mathbb Z/2) = \mathbb Z/2$. Thus exactly one of $\beta^*$ and $\beta^* + \beta$ is zero in $H_1(M; \mathbb Z/2)$.
(Recall that duality implies that the image of $H_1(\partial M; \mathbb Z/2)$
in $H_1(M; \mathbb Z/2)$ is $\mathbb Z/2$.) Since $\beta$ lies in the
kernel of $\rho_2$, it follows that $\rho_2(\beta^*)$ is sent to
zero in $H_1(D_3; \mathbb Z/2)$.
But then $\rho_2(\beta^*)$ has order $3$ in $D_3$, contrary to what we deduced in the previous paragraph. Thus $s_0 = 2$. Now apply the argument at the end of the proof of Proposition 6.6 of \cite{BCSZ2} to see that $4 = 2s_0 \geq |\partial P| \geq 4$. Hence $P$ has four boundary components.
\qed

The four-punctured $2$-sphere $P$ cuts $M$ into two components $X_1$ and $X_2$. If $P_i$ denotes the copy of $P$ in $\partial X_i$ then $M$ is the union of $X_1$ and $X_2$ with $P_1$ and $P_2$ identified by a homeomorphism $f: P_1\to P_2$. The boundary of $P$ cuts $\partial M$ into four annuli $A_{11}, A_{21}, A_{12}, A_{22}$ listed in the order they appear around $\partial M$,
where $A_{11}, A_{12}$ are contained in $X_1$ and
$A_{21}, A_{22}$ are contained in $X_2$.
The arguments given in the proof of Lemma 4.5 of \cite{BCSZ2} show that for
each $i$, the two annuli $A_{i1}$ and $A_{i2}$ in $X_i$ are unknotted
and unlinked.
This means that there is a neighbourhood of $A_{i1}\cup A_{i2}$ in $X_i$
which is homeomorphic to $E_i\times I$, where $E_i$ is a thrice-punctured $2$-sphere and $I$ is the interval $[0,1]$, such that $(E_i \times  I) \cap P_i = (E_i \times \partial I)$, and the exterior of $E_i\times I$ in $X_i$ is a solid torus $V_i$. We label  the boundary components of $E_i$ as $\partial_j E_i$ ($j=1,2,3$) so that $\partial_j E_i \times I = A_{ij}$ for $j=1,2$.

Let $\hat P$ be the two sphere in $M(\beta)$ obtained from $P$ by  capping off $\partial P$
with four meridian  disks from the filling solid torus $V_\beta$.
These disks cut $V_\beta$ into four $2$-handles $H_{ij}$ ($i,j =1,2$) such that the attaching annulus  of $H_{ij}$ is $A_{ij}$ for each $i,j$.
Let $X_i(\beta)$ be the manifold obtained by attaching $H_{ij}$ to $X_i$ along
$A_{ij}$ ($j=1,2$). Then $X_1(\beta)$ is a once-punctured $L(2,1)$
and $X_2(\beta)$ a once-punctured  $L(3,1)$.

It follows from the description above that $X_1$ is obtained from
$E_1 \times I$ and $V_1$ by identifying $\partial_3 E_1 \times I$ with
an annulus $A_1$ in $\partial V_1$ whose core curve is a $(2,1)$ curve
in $\partial V_1$.
We can assume that $A_1$ is invariant under the standard involution of $V_1$ whose fixed point set is a pair of arcs contained in disjoint meridian disks of $V_1$. Note that the two boundary components of $A_1$ are interchanged under this map. Similarly, $X_2$ is obtained from $E_2 \times I$ and $V_2$ by identifying $\partial_3 E_2 \times I$ with an annulus $A_2$ in $\partial V_2$ whose core curve is a $(3,1)$ curve in $\partial V_2$. Again we can suppose that $A_2$ is invariant under the standard involution of $V_2$ which interchanges the two boundary components of $A_2$. See Figures \ref{inv1} and \ref{inv2}.

The map $f|: \partial P_1 \to \partial P_2$ is constrained in several ways by our hypotheses. For instance, the fact that $\partial M$ is connected implies that $f(\partial_1 E_1 \times \{i\}) = \partial_2 E_2 \times \{j\}$ for some $i, j$. Other conditions are imposed by the homology of $M$.

\begin{lemma} \label{gluing}
We can assume that either \\
$(a)$ $f(\partial_1 E_1 \times \{0\}) = \partial_1 E_2 \times \{0\}, f(\partial_2 E_1 \times \{0\}) = \partial_2 E_2 \times \{0\}, f(\partial_1 E_1 \times \{1\}) = \partial_2 E_2 \times \{1\}$,  and $f(\partial_2 E_1 \times \{1\}) = \partial_1 E_2 \times \{1\}$, or \\
$(b)$  $f(\partial_1 E_1 \times \{0\}) = \partial_1 E_2 \times \{0\}, f(\partial_2 E_1 \times \{0\}) = \partial_1 E_2 \times \{1\}, f(\partial_1 E_1 \times \{1\}) = \partial_2 E_2 \times \{1\},$ and $f(\partial_2 E_1 \times \{1\}) = \partial_2 E_2 \times \{0\}$.
\end{lemma}

\pf Without loss of generality we can suppose that $f(\partial_1 E_1 \times \{0\}) = \partial_1 E_2 \times \{0\}$. Hence, as $\partial M$ is connected, one of the following four possibilities arises:
\vspace{-.35cm}
\begin{itemize}
\item[(a)] $f(\partial_2 E_1 \times \{0\}) = \partial_2 E_2 \times \{0\}, f(\partial_1 E_1 \times \{1\}) = \partial_2 E_2 \times \{1\},$ and $f(\partial_2 E_1 \times \{1\}) = \partial_1 E_2 \times \{1\}$.
\item[(b)] $f(\partial_2 E_1 \times \{0\}) = \partial_1 E_2 \times \{1\}, f(\partial_1 E_1 \times \{1\}) = \partial_2 E_2 \times \{1\},$ and $f(\partial_2 E_1 \times \{1\}) = \partial_2 E_2 \times \{0\}$.
\item[(c)] $f(\partial_2 E_1 \times \{0\}) = \partial_1 E_2 \times \{1\}, f(\partial_1 E_1 \times \{1\}) = \partial_2 E_2 \times \{0\},$ and $f(\partial_2 E_1 \times \{1\}) = \partial_2 E_2 \times \{1\}$, or
\item[(d)] $f(\partial_2 E_1 \times \{0\}) = \partial_2 E_2 \times \{1\}, f(\partial_1 E_1 \times \{1\}) = \partial_2 E_2 \times \{0\},$ and $f(\partial_2 E_1 \times \{1\}) = \partial_1 E_2 \times \{1\}$.
\end{itemize}
\vspace{-.35cm}


Let $a_i,b_i ,x_i \in H_1(X_i)$ be represented, respectively,
by $\partial_1E_i$, $\partial_2E_i$, and a core of $V_i$, $i=1,2$.
Then $H_1(X_i)$ is the abelian group generated  by $a_i,b_i,x_i$,
subject to the relation
\begin{gather}
2x_1 = a_1 + b_1\ ,\qquad i=1 \tag{i}\\
3x_2 = a_2 + b_2\ ,\qquad i=2 \tag{ii}
\end{gather}
Since $f(\partial_1 E_1\times \{0\})  = \partial_1 E_2 \times\{0\}$,
we may orient $\partial E_1$, $\partial E_2$ so that in $H_1 (M)$ we
have $a_1=a_2$.
Then $H_1 (M)$ is the quotient of $H_1 (X_1) \oplus H_1 (X_2)$ by this
relation together with the additional relations corresponding to the
four possible gluings:
\begin{itemize}
\item[(a)] $b_1 = b_2$, $b_1 = a_2$
\item[(b)] $b_1 =  - a_2$, $b_1 = -b_2$
\item[(c)] $b_1 = - a_2$, $b_1 = b_2$
\item[(d)] $b_1 = - b_2$, $b_1 = a_2$
\end{itemize}
Taking $\zed/3$ coefficients, equation (i) allows us to eliminate
$x_1$, while (ii) gives $a_2 + b_2 =0$.
Hence $H_1 (M;\zed/3) \cong\zed/3 \oplus A$, where the $\zed/3$
summand is generated by $x_2$ and $A$ is defined by generators
$b_1,a_2,b_2$, and relations $a_2 + b_2 =0$ plus those listed in
(a), (b), (c) and (d) above.
Thus $A=0$ in cases~(a) and (b), and $A\cong\zed/3$ in cases~(c) and (d).
Since $H_1 (M) \cong \zed\oplus \zed/2$, we conclude that cases~(c) and
(d) are impossible.\qed

\section{The proof of Theorem \ref{rf} when case (a) of Lemma \ref{gluing} arises}
\label{case a}

\begin{figure}[!ht]
{\epsfxsize=4in \centerline{\epsfbox{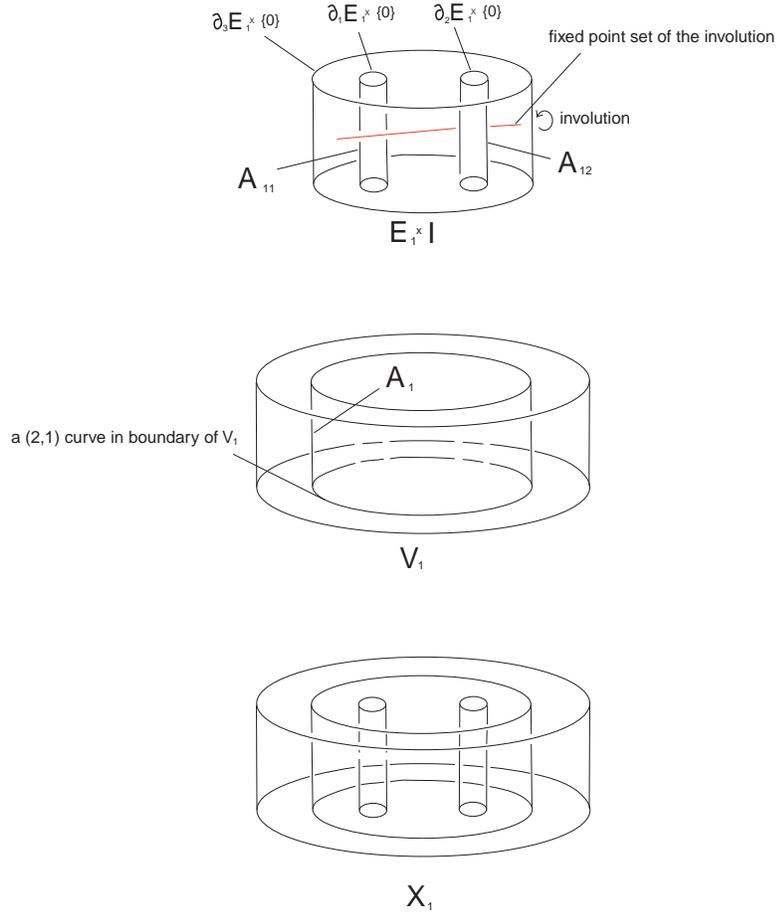}}\hspace{10mm}}
\caption{$X_1$ and the involution $\tau_1$}\label{inv1}
\end{figure}

Figure \ref{inv1} depicts an involution $\tau_1$ on $E_1 \times I$ under which $\partial_3 E_1\times I$ is invariant, has its boundary components interchanged, and $\tau_1(A_{11})=A_{12}$. Then $\tau_1$ extends to an involution of  $X_1$ since its restriction to $\partial_3 E_1\times I = A_1$ coincides with the restriction to $A_1$ of the standard involution of $V_1$. Evidently $\tau_1(\partial_1 E_1 \times \{0\}) = \partial_2 E_1 \times \{1\}$ and $\tau_1(\partial_2 E_1 \times \{0\}) = \partial_1 E_1 \times \{1\}$.

Figure \ref{inv2} depicts an involution $\tau_2$ on $E_2 \times I$ under
which each of the annuli $\partial_3 E_2 \times I, A_{21}$, and $A_{22}$
are invariant.
Further, it interchanges the components of $E_2 \times\partial I$ and as in
the previous paragraph, $\tau_2$ extends to an involution of  $X_2$.
Note that $\tau_2(\partial_j E_2 \times \{0\}) = \partial_j E_2 \times \{1\}$
for $j = 1, 2$.

Next consider the orientation preserving involution $\tau_2' = f (\tau_1|P_1) f^{-1}$ on $P_2$. By construction we have $\tau_2'(\partial_j E_2 \times \{0\}) = \partial_j E_2 \times \{1\}$ for $j = 1, 2$, and therefore $\tau_2' = g (\tau_2|P_2) g^{-1}$ where $g: P_2 \to P_2$ is a homeomorphism whose restriction to $\partial P_2$ is isotopic to $1_{\partial P_2}$. The latter fact implies that $g$ is isotopic to a homeomorphism $g': P_2 \to P_2$ which commutes with $\tau_2|P_2$. Hence, $\tau_2'$ is isotopic to $\tau_2|P_2$ through orientation preserving involutions whose fixed point sets consist of two points. In particular, $\tau_1$ and $\tau_2$ can be pieced together to form an orientation preserving involution $\tau: M \to M$.

For each slope $\gamma$ on $\partial M$, $\tau$ extends to an involution $\tau_\gamma$ of the associated Dehn filling $M(\gamma) = M \cup V_\gamma$, where $V_\gamma$ is the filling solid torus. Thurston's orbifold theorem applies to our situation and implies that $M(\gamma)$ has a geometric decomposition. In particular, $M(\alpha)$ is a Seifert fibred manifold whose base orbifold is of the form $S^2(2,3,4)$, a $2$-sphere with three cone points of orders $2,3,4$ respectively.

It follows immediately from our constructions that $X_1(\beta)/ \tau_\beta$
and $X_2(\beta)/ \tau_\beta$ are $3$-balls.
Thus $M(\beta)/ \tau_\beta = (X_1(\beta)/ \tau_\beta) \cup (X_2(\beta)/
\tau_\beta) \cong S^3$ and since $\partial M / \tau \cong S^2$, it follows
that $M / \tau$ is a $3$-ball.
More precisely, $M/\tau$ is an orbifold $(N,L^0)$, where $N$ is a 3-ball,
$L^0$ is a properly embedded 1-manifold in $N$ that meets $\partial N$ in
four points, and $M$ is the double branched cover of $(N,L^0)$.
We will call $(N,L^0)$ a {\em tangle}, and if we choose some identification
of $(\partial N,\partial L^0)$ with a standard model of ($S^2$, {\em four
points}), then $(N,L^0)$ becomes a {\em marked} tangle.
Capping off $\partial N$ with a 3-ball $B$ gives $N\cup_\partial B\cong S^3$.
Then, if $\gamma$ is a slope on $\partial M$, we have $V_\gamma/\tau_\gamma
\cong (B,T_\gamma)$, where $T_\gamma$ is the rational tangle in $B$
corresponding to the slope $\gamma$.
Hence
\begin{equation*}
\begin{split}
M(\gamma)/\tau_\gamma & = (M/\tau) \cup (V_\gamma/\tau_\gamma)\\
& = (N,L^0) \cup (B,T_\gamma)\\
& = (S^3, L^0 (\gamma))\ ,
\end{split}
\end{equation*}
where $L^0(\gamma)$ is the link in $S^3$ obtained by capping off $L^0$ with
the rational tangle $T_\gamma$.

\begin{figure}[!ht]
{\epsfxsize=4in \centerline{\epsfbox{inv2.ai}}\hspace{10mm}}
\caption{$X_2$ and the involution $\tau_2$}\label{inv2}
\end{figure}

\begin{figure}[!ht]
{\epsfxsize=5in \centerline{\epsfbox{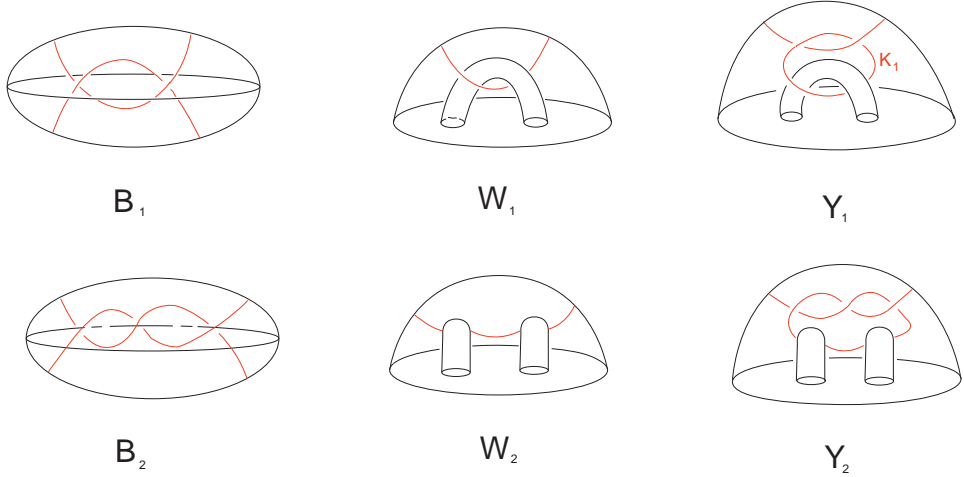}}\hspace{10mm}}
\caption{The branch sets in $B_i, W_i$, and $Y_i$}\label{quo1}
\end{figure}

\begin{figure}[!ht]
{\epsfxsize=5in \centerline{\epsfbox{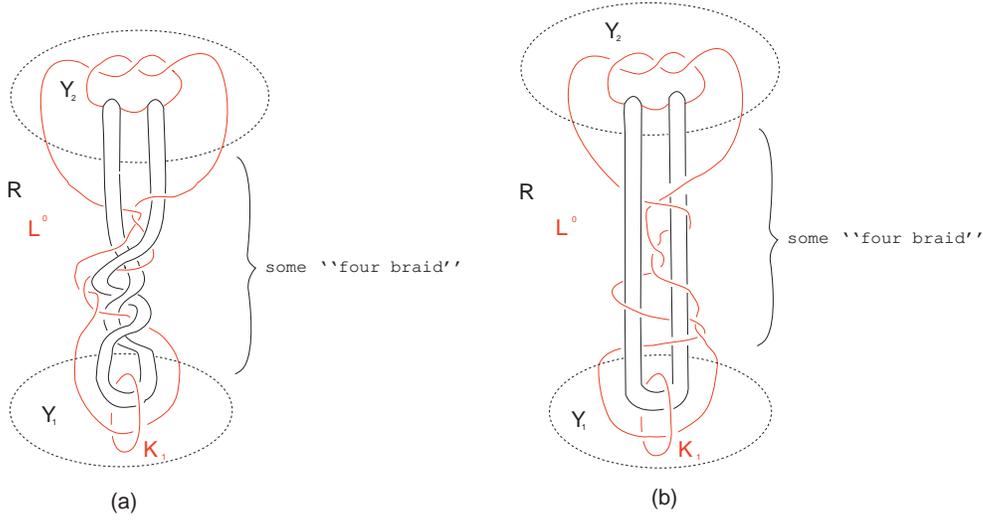}}\hspace{10mm}}
\caption{The branch set $L^0$ in $N$}\label{filling}
\end{figure}

We now give a more detailed description of the tangle $(N,L^0)$.
For $i=1,2$, let $B_i=V_i/ \tau_i$, $W_i=E_i\times I/\tau_i$,
$Y_i=X_i/\tau_i$, and  $Q_i=P_i/\tau_i$.
Figure \ref{quo1}  gives a detailed description of the branch sets
 in $B_i$, $W_i$,  $Y_i$ with respect to the
corresponding branched covering maps.
Note that $N$ is the union of
$Y_1, Y_2$, and a product region $R \cong Q_1 \times I$ from $Q_1$ to $Q_2$
which intersects the branch set $L^0$ of the cover $M \to N$ in a $2$-braid.
In fact, it is clear from our constructions that we can think of the union
$(L^0 \cap R) \cup (\partial N \cap R)$ as a ``$4$-braid" in $R$ with two
``fat strands" formed by $\partial N \cap R$.
See Figure \ref{filling}(a).
By an isotopy of $R$ fixing $Q_2$, and which keeps $R$, $Q_1$, and $Y_1$
invariant, we may untwist the crossings between the two fat strands in
Figure \ref{filling}(a) so that the pair $(N, L^0)$ is as depicted in
Figure \ref{filling}(b).

The slope $\beta$ is the boundary slope of the planar surface $P$, and
hence the rational tangle $T_\beta$ appears in Figure~\ref{filling}(b)
as two short horizontal arcs in $B$ lying entirely in $Y_2(\beta) =
X_2 (\beta)/\tau_\beta$.
Since $\Delta (\alpha,\beta)=2$, $T_\alpha$ is a tangle of the form shown
in Figure~\ref{filling6}(a).
Recall that $M(\alpha)$ is a Seifert fibred manifold with base orbifold of
type $S^2 (2,3,4)$, and is the double branched cover of $(S^3,L^0(\alpha))$.
Write $L= L^0(\alpha)$.

\begin{lemma} \label{type}
$L$ is a Montesinos link of type $(\frac{p}{2}, \frac{q}{3}, \frac{r}{4})$.
\end{lemma}

\pf By Thurston's orbifold theorem, the Seifert fibering of $M(\alpha)$ can be isotoped to be
invariant under $\tau_\alpha$. Hence the quotient orbifold is Seifert fibered
in the sense of Bonahon-Siebenmann, and so either $L$ is a Montesinos link
or $S^3 \setminus L$ is Seifert fibred. From Figure 6(a) we see that $L$ is
a 2-component link with an unknotted component and linking number $\pm 1$.
But the only link $L$ with this property such that $S^3 \setminus L$ is
Seifert fibred is the Hopf link (see \cite{BM}), whose $2$-fold cover is
$P^3$. Thus $L$ must be a Montesinos link. Since the base orbifold of
$M(\alpha)$ is $S^2(2,3,4)$, $L$ has type $(\frac{p}{2}, \frac{q}{3},
\frac{r}{4})$ (c.f \S 12.D of \cite{BuZi}).\qed

It's easy to check that any Montesinos link $L$ of the type described in the
Lemma \ref{type} has two components, one of which, say $K_1$, is a trivial
knot, and the other, $K_2$, a trefoil knot.
Our goal is to use the particular nature of our situation to show that the
branch set $L$ cannot be a Montesinos link of type
$(\frac{p}{2}, \frac{q}{3}, \frac{r}{4})$, and thus derive a contradiction.


{From} Figure \ref{filling}, we see that $L^0$ has a closed, unknotted
component, which must be the component $K_1$  
of the Montesinos link of type $(\frac{p}{2}, \frac{q}{3}, \frac{r}{4})$
described above.
Then $L^0 \setminus K_1 = K_2 \cap N$, which we denote by $K_2^0$.

Now  delete $K_1$ from $N$ and let $U$ be the double branched cover
of $N$ branched over $K_2^0$. Then $U$ is a compact, connected,
orientable $3$-manifold with boundary a torus which can be
identified with $\partial M$. In particular, if we consider $\alpha$
and $\beta$ as slopes on $\partial U$, then both $U(\alpha)$ and
$U(\beta)$ are the lens space $L(3,1)$, since they are $2$-fold
covers of $S^3$ branched over a trefoil knot. Hence the cyclic
surgery theorem of \cite{CGLS} implies that $U$ is either a  Seifert
fibred space or a reducible manifold.

\begin{lemma}
$U$ is not a Seifert fibred space.
\end{lemma}

\pf Suppose $U$ is a Seifert fibred space, with base surface $F$ and $n \geq 0$
exceptional fibres. If $F$ is non-orientable then $U$ contains a Klein bottle,
hence $U(\alpha) \cong L(3,1)$ does also. But since non-orientable surfaces in $L(3,1)$ are non-separating, this implies that $H_1(L(3,1); \mathbb Z/2) \not \cong 0$, which is clearly false. Thus $F$ is orientable.

If $U$ is a solid torus then clearly $U(\alpha) \cong U(\beta) \cong L(3,1)$ implies
$\Delta(\alpha, \beta) \equiv 0$ (mod $3$), contradicting the fact that $\Delta(\alpha,
\beta) = 2$. Thus we assume that $U$ is not a solid torus, and take $\phi \in H_1(\partial U)$ to be the slope on $\partial U$ of a Seifert fibre. Then $U(\phi)$ is reducible \cite{Hl} so $d = \Delta(\alpha, \phi) > 0$, and $U(\alpha)$ is a Seifert fibred space with base surface $F$ capped off with a disk, and $n$ or $n+1$ exceptional fibres, according as $d = 1$ or $d > 1$. Since $U(\alpha)$ is a lens space and $U$ isn't a solid torus, we must have that $F$ is a disk, $n = 2$, and $d = 1$. Similarly $\Delta (\beta, \phi) = 1$. In particular, without loss of generality $\beta = \alpha + 2\phi$ in $H_1(\partial U)$.

The base orbifold of $U$ is of the form $D^2(p,q)$, with $p,q > 1$. Then $H_1(U)$ is the abelian group defined by generators $x,y$ and the single relation $px + qy = 0$. Suppose $\alpha \mapsto ax + by$ in $H_1(U)$. Then $H_1(U(\alpha))$ is presented by the matrix $\left( \begin{array}{cc} p & a \\ q & b \end{array} \right)$. Similarly, since $\phi \mapsto px$ in $H_1(U)$, $H_1(U(\beta))$ is presented by $\left( \begin{array}{cc} p & a + 2p\\ q & b \end{array} \right)$. But the determinants of these matrices differ by $2pq \geq 8$, so they cannot both be $3$ in absolute value. This completes the proof of the lemma.
\qed

Thus $U$ is reducible, say $U \cong V \# W$ where $\partial V = \partial U$ and $W \not \cong S^3$ is closed. Consideration of $M(\alpha)$ and $M(\beta)$ shows that $W \cong L(3,1)$ and $V(\alpha) \cong V(\beta) \cong S^3$, and so Theorem 2 of \cite{GL1} implies that $V \cong S^1 \times D^2$. It follows that any simple closed curve in $\partial U$ which represents either $\alpha$
or $\beta$ is isotopic to the core curve of $V$.
Let $\lambda \in H_1(\partial U)$ denote the
meridional slope of $V$.
Then $\{\beta, \lambda\}$ is a basis of $H_1(\partial U)$ and up to changing
the sign of $\alpha$ we have $\alpha=\beta \pm 2\lambda$.
%
%

Since $U\cong (S^1 \times D^2)\,\#\, L(3,1)$, we can find a
homeomorphism between the pair $(N,K_2^0)$ and the tangle shown in
Figure~\ref{NK20}(a), with the $\beta$, $\alpha$, and $\lambda$
fillings shown in Figures~\ref{NK20}(b), (c) and (d) respectively.
(We show the case $\alpha = \beta +2\lambda$; the other possibility
can be handled similarly.)
\begin{figure}[!ht]
{\epsfxsize=4in \centerline{\epsfbox{NK20.ai}}\hspace{10mm}}
\caption{}      
\label{NK20}
\end{figure}


Recall that in Figure \ref{filling}(b), the slope $\beta$
corresponds to the rational tangle consisting of two short
``horizontal'' arcs in the filling ball $B$. It follows that under
the homeomorphism from the tangle shown in Figure~\ref{NK20}(a) to
$(N,K_2^0)$ shown in Figure~\ref{filling}(b),
the tangle $T_\alpha$ (resp. $T_\lambda$) is sent
to a rational tangle of the form
shown in Figure~\ref{filling6}(a) (resp. \ref{filling6}(b)).
{From} Figure~\ref{NK20}(d) we see that $L^0(\gamma)$ is
a link of three components $K_1 \cup O_1 \cup K_3$, where $O_1$ is a
trivial knot which bounds a disk $D$ disjoint from $K_3$ and which
intersects $\partial N$ in a single arc; see
Figure~\ref{filling6}(b). Push the arc $O_1\cap B$ with its two
endpoints fixed into $\p B$ along $D$, and let $O_1^*$ be the
resulting  knot (see part (c) of Figure \ref{filling6}). Then there
is a disk $D_*$ (which is a subdisk of $D$) satisfying the following
conditions:
\begin{itemize}
\item[(1)] $\partial D_*=O_1^*$.
\item[(2)] $D_*$ is  disjoint from $K_3$.
\item[(3)] The interior of $D_*$ is disjoint from $B$.
\end{itemize}
\begin{figure}[!ht]
{\epsfxsize=6in \centerline{\epsfbox{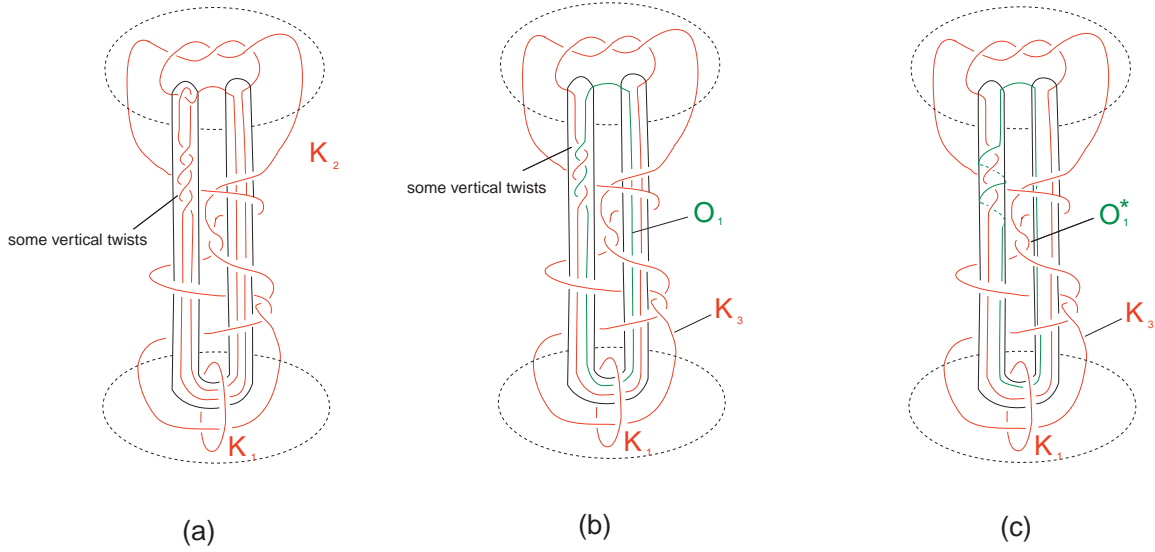}}\hspace{10mm}}
\caption{The tangle fillings $N(\alpha)$  and $N(\lambda)$.}
\label{filling6}
\end{figure}

\noindent Perusal of Figure \ref{filling6} (c) shows that the following condition is
also achievable.
\begin{itemize}
\item[(4)] $D_* \cap Q_2$ has a single arc component, and this arc
component connects the two boundary components of $Q_2$ and is outermost
in $D_*$ amongst the components of  $D_*  \cap Q_2$.
\end{itemize}
Among all disks in $S^3$ which satisfy conditions
(1)--(4), we may assume that $D_*$ has been chosen so that
\begin{itemize}
\item[(5)] $D_* \cap Q_2$ has the minimal number of components.
\end{itemize}

\begin{claim}\label{parallel}
Suppose that $D_* \cap Q_2$ has  circle components. Then each such
circle separates $K_3 \cap Q_2$ from $\partial Q_2$ in $Q_2$.
\end{claim}

\pf Let $\delta$ be a circle component of $D_* \cap Q_2$.
Then $\delta$ is essential in $Q_2 \setminus (Q_2 \cap K_3)$, for if it
bounds a disk $D_0$ in $Q_2 \setminus (Q_2 \cap K_3)$, then an innermost
component of $D_* \cap D_0 \subset D_* \cap Q_2$ will bound a disk
$D_1 \subset D_0$.
We can surger $D_*$ using $D_1$ to get a new disk satisfying conditions
(1)--(4)  above, but with  fewer components of intersection
with $Q_2$ than $D_*$, contrary to assumption (5).

Next since the arc component of $D_* \cap Q_2$ connects the two
boundary components of $Q_2$, $\delta$ cannot separate the two
boundary components of $Q_2$ from each other.

Lastly suppose that $\delta$ separates the two points of $Q_2 \cap K_3$.
Then $\delta$ is isotopic to a meridian curve of $K_3$ in $S^3$.
But this is impossible since $\delta$ also bounds a disk in $D_*$ and is
therefore null-homologous in $S^3 \setminus K_3$.
The claim follows.\qed
\begin{figure}[!ht]
{\epsfxsize=4in \centerline{\epsfbox{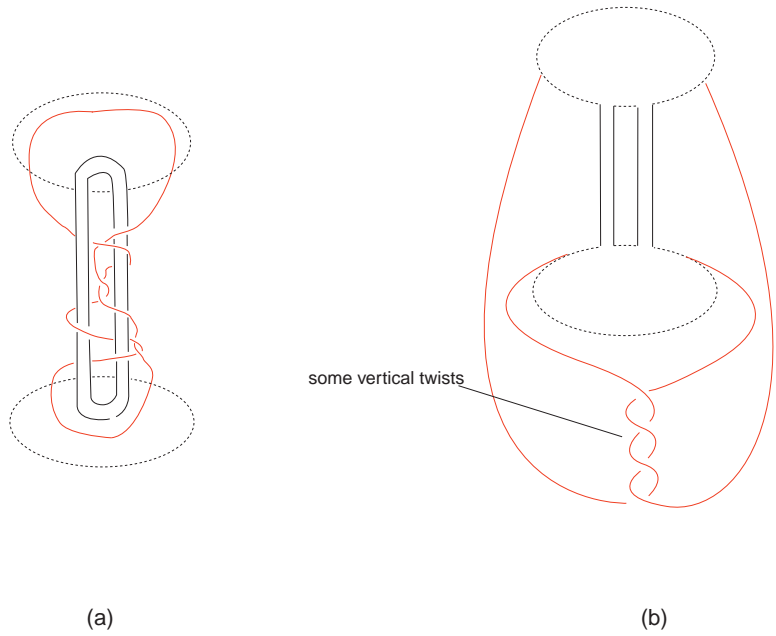}}\hspace{10mm}}
\caption{Capping off the $4$-braid to obtain a trivial
link}\label{filling7}
\end{figure}

It follows from Claim \ref{parallel} that there are disjoint arcs in
$Q_2$, one, say $\sigma_1$, which connects the two points of $Q_2
\cap K_3$ and is disjoint from $D_*$, and $\sigma_2 = D_* \cap Q_2$
the other. Hence we obtain a ``2-bridge link'' of two components --
one fat, one thin -- in $S^3$ by capping off the ``$4$-braid'' in
$R$ with $\sigma_1$ and $\sigma_2$ in $Y_2 \subset Y_2(\beta)$ and
with $K_3 \cap Y_1 \subset Y_1$ and $\partial N \cap Y_1 \subset
Y_1$ in the $3$-ball $Y_1(\beta)$ (see Figure \ref{filling7}(a)).
Furthermore, since the disk $D_*$  gives a disk bounded by the ``fat
knot'' which is disjoint from the ``thin knot'', the link is a
trivial link.
\begin{figure}[!ht]
{\epsfxsize=4in \centerline{\epsfbox{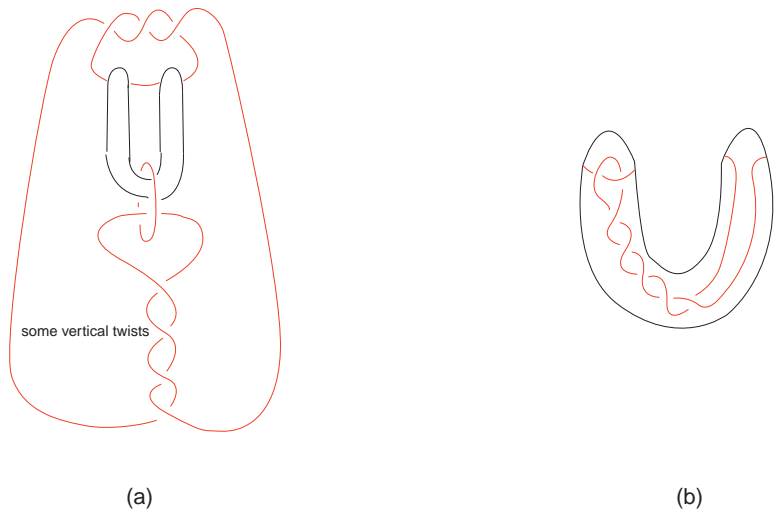}}\hspace{10mm}}
\caption{The pair $(N, L^0)$ and the filling tangle
$T_\a$}\label{filling8}
\end{figure}
\begin{figure}[!ht]
{\epsfxsize=4in \centerline{\epsfbox{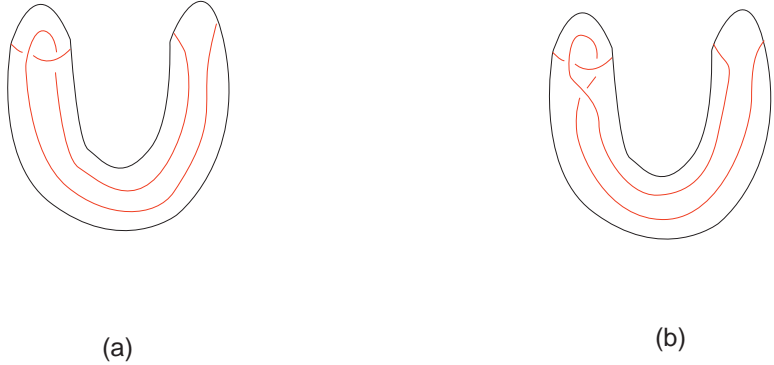}}\hspace{10mm}}
\caption{The two possible $T_\a$}\label{filling8b}
\end{figure}

\begin{figure}[!ht]
{\epsfxsize=5in \centerline{\epsfbox{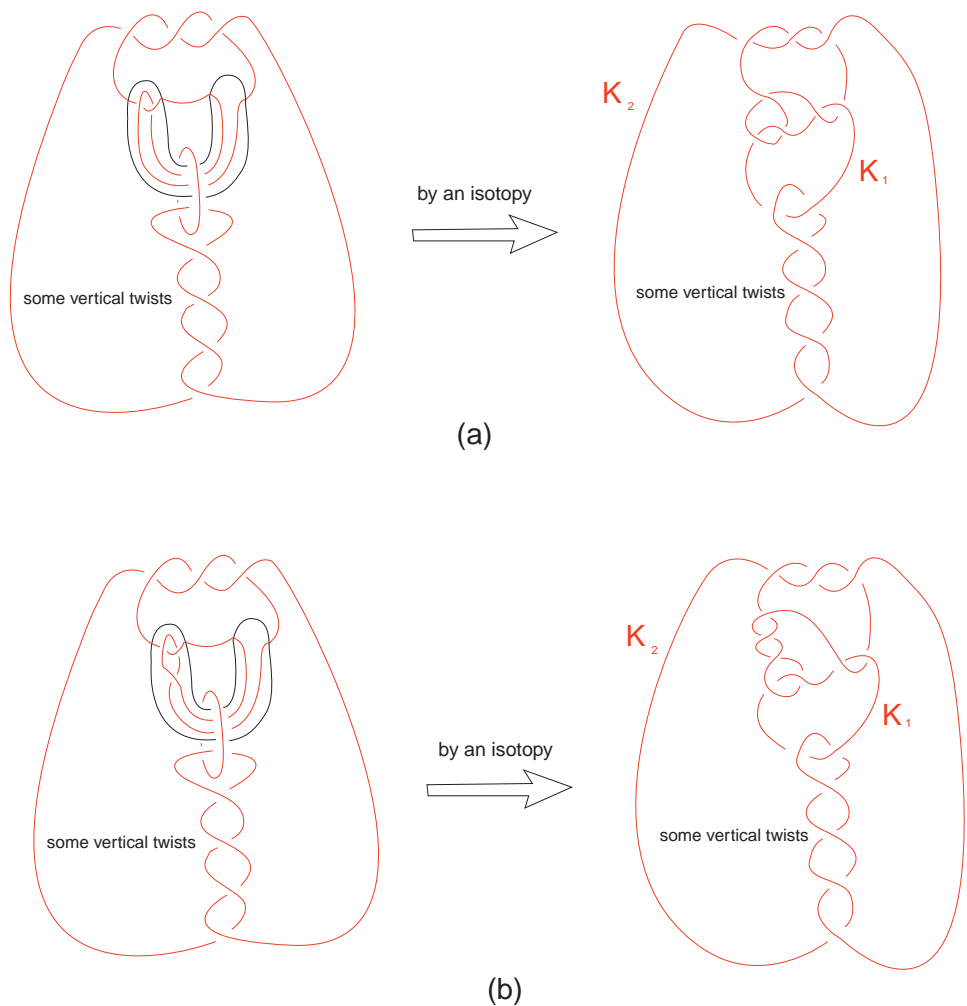}}\hspace{10mm}}
\caption{$L$ as a Montesinos link of the type $(\frac{1}{3},
\frac{-3}{8}, \frac{m}{2})$ or $(\frac13, \frac{-5}{8},
\frac{m}{2})$}\label{filling9}
\end{figure}

Now it follows from the standard presentation of a $2$-bridge link
as a $4$-plat (see \S 12.B of \cite{BuZi}), that there is an isotopy
of $R$, fixed on the ends $Q_1,Q_2$ and on the two fat strands,
taking the ``4-braid'' to one of the form shown in
Figure~\ref{filling7}(b). Hence $(N,L^0)$ has the form shown in
Figure~\ref{filling8}(a). The filling rational tangle $T_\a$ is of
the form  shown in Figure~\ref{filling8}(b). Since the component
$K_2^0 (\alpha)$ of $L^0(\alpha) =L$ has to be a trefoil, there are
only two possibilities for the number of twists in $T_\alpha$; see
Figure \ref{filling8b}. The two corresponding possibilities for $L$
are shown in Figure~\ref{filling9}. But these are Montesinos links
of the form $(\frac13, \frac{-3}8,\frac{m}2)$ and $(\frac13,
\frac{-5}8, \frac{m}2)$, respectively.

This final contradiction completes the proof of Theorem~\ref{rf} under the assumptions of case (a) of Lemma \ref{gluing}.\qed

\section{The proof of Theorem \ref{rf} when case (b) of
Lemma \ref{gluing} arises}

\begin{figure}[!ht]
{\epsfxsize=5.5in \centerline{\epsfbox{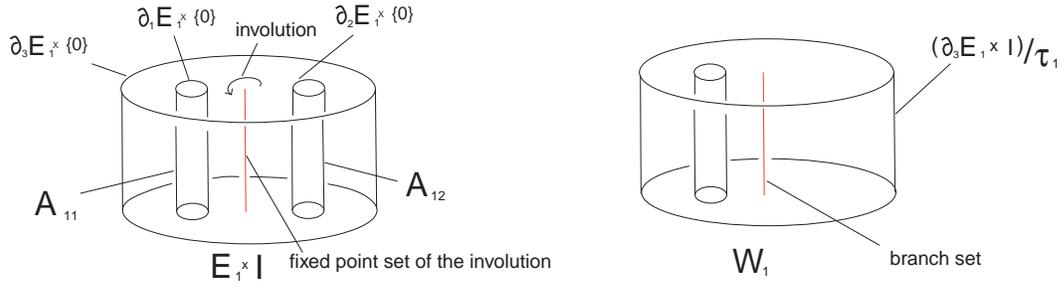}}\hspace{10mm}}
\caption{Involution on $E_1\times I$}\label{inv3}
\end{figure}

In this case we choose an involution $\tau_1$ on $E_1 \times I$  as
shown in Figure \ref{inv3}.
Then $\tau_1(\partial_3 E_1 \times \{j\}) = \partial_3 E_1 \times \{j\}$,
$\tau_1(\partial_1 E_1\times \{j\}) = \partial_2 E_1\times \{j\}$ ($j = 0,1$),
and the restriction of $\tau_1$ on
$\partial_3 E_1 \times I$ extends to an involution of $V_1$
whose fixed point set is a core circle of this solid torus.
Thus we obtain an
involution $\tau_1$ on $X_1$. The quotient of $V_1$ by $\tau_1$ is a
solid torus $B_1$ whose core circle is the branch set.
Further, $A_1 / \tau_1$ is a longitudinal annulus of $B_1$.
The quotient of $E_1 \times I$ by $\tau_1$ is also solid torus $W_1$ in
which $(\partial_3 E_1 \times I) / \tau_1$ is a longitudinal annulus.
Figure \ref{inv3} depicts $W_1$ and its branch set.
It follows that the pair $(Y_1 = X_1/ \tau_1, \hbox{branch set of } \tau_1)$
is identical to the analogous pair in Section~\ref{case a}
(see Figure~\ref{quo1}).

Next we take $\tau_2$ to be the same involution on $X_2$ as that
used in Section~\ref{case a}. An argument similar to the one used in
that section shows that $\tau_1$ and $\tau_2$ can be pieced together
to form an involution $\tau$ on $M$. {From} the previous paragraph
we see that the quotient $N = M / \tau$ and its branch set are the
same as those in Section~\ref{case a}. Hence the argument of that
section can be used from here on to obtain a contradiction. This
completes  the proof of Theorem \ref{rf} in case (b). \qed

\newpage

\def\bysame{$\underline{\hskip.5truein}$}

{\small
\vspace{5mm}
\noindent
D\'epartement de math\'ematiques, Universit\'e du Qu\'ebec \`a Montr\'eal,
P. O. Box 8888, Centre-ville, Montr\'eal, Canada, H3C 3P8 \\
e-mail: steven.boyer@uqam.ca

\vspace{2mm}
\noindent
Department of Mathematics, University of Texas at Austin,
Austin, TX, 78712-1082, USA \\
e-mail: gordon@math.utexas.edu

\vspace{2mm}
\noindent
Department of Mathematics, University at Buffalo,
Buffalo, NY 14260-2900, USA  \\
e-mail: xinzhang@buffalo.edu}
\end{document}